\documentclass[11pt, reqno]{amsart}
\usepackage{natbib}

\theoremstyle{definition}
\newtheorem{definition}{Definition}[section]
\newtheorem{example}[definition]{Example}
\newtheorem{remark}[definition]{Remark}
\theoremstyle{plain}
\newtheorem{theorem}[definition]{Theorem}
\newtheorem{lemma}[definition]{Lemma}

\newtheorem{proposition}[definition]{Proposition}

\begin{document}
\title[Right Division in Groups]{Right Division in Groups,
Dedekind-Frobenius Group Matrices, and Ward Quasigroups}

\author{Kenneth W. Johnson}

\address{Penn State Abington, 1600 Woodland Rd, Abington, PA 19001, USA}

\email{kwj1@psu.edu}

\author{Petr Vojt\v{e}chovsk\'y}

\address{Department of Mathematics, University of Denver, 2360 S Gaylord St,
Denver, CO 80208, USA}

\email{petr@math.du.edu}



\thanks{The second author partially supported by Grant Agency of Charles
University, grant number $269/2001/$B-MAT/MFF}

\begin{abstract}
The variety of quasigroups satisfying the identity $(xy)(zy)=xz$ mirrors the
variety of groups, and offers a new look at groups and their multiplication
tables. Such quasigroups are constructed from a group using right division
instead of multiplication. Their multiplication tables consist of circulant
blocks which have additional symmetries and have a concise presentation.

These tables are a reincarnation of the group matrices which Frobenius used
to give the first account of group representation theory. Our results imply
that every group matrix may be written as a block circulant matrix and that
this result leads to partial diagonalization of group
matrices, which are present in modern applied mathematics. We also discuss
right division in loops with the antiautomorphic inverse property.
\end{abstract}

\maketitle



\section{Introduction}

\noindent Several papers have characterized groups using the operation of right division
 $x\cdot y^{-1}$ instead of the multiplication $x\cdot y$.
 However, it
is not clear how much is gained within group theory per se by such a change in perspective. The aim of
this paper is to suggest that there are some advantages, connected to the added
symmetry of the multiplication table.

Right division was used already by Frobenius. In his first papers on group
representation theory \cite{F1}, \cite{F2}, the essential objects were the
group matrix and its determinant, the group determinant. For a finite group
$G=\{g_{1}$, $\dots $, $g_{n}\}$, its \emph{group matrix} $X_{G}$ is defined to
be the $n\times n$ matrix whose $(i,j)$th entry is $x_{g_{i}g_{j}^{-1}}$, where
$\{x_{g_{1}}$, $\dots $, $x_{g_{n}}\}$ is a set of commuting variables. Usually
the term group matrix is also applied to any matrix obtained from a group
matrix by assigning values in a ring to the variables. We refer to \cite{J1},
\cite{J2} and \cite{HJ1} for information on how Frobenius' ideas have
stimulated recent research. Frobenius relied heavily on the symmetrical nature
of the group matrix in his proofs of the basic results of representation
theory.

Examples of group matrices were known well before Frobenius, since a
\emph{circulant matrix}, i.e., a matrix of the form
\begin{equation*}
\left(
\begin{array}{cccc}
c_1 & c_2 & \cdots & c_n \\
c_n & c_1 & \cdots & c_{n-1} \\
\vdots & \vdots &  & \vdots \\
c_2 & c_3 & \dots & c_1
\end{array}
\right),
\end{equation*}
is a group matrix of the cyclic group $C_n$ of order $n$. We denote this
circulant by $C(c_1,\dots,c_n)$. The efficient diagonalization of circulants is
behind the Finite Fourier Transform \cite{AT}. We refer to \cite{dav} for a
thorough account of circulants. Group matrices of arbitrary groups have
appeared in probability \cite{dia}, where they arise as transition matrices for
Markov chains. They have also appeared in the theory of tight frames as Gram
matrices \cite{wal} with connections to wavelets and non-harmonic Fourier
series.

We may interpret the matrix $X_{G}$ as an encoding of the multiplication table
of a quasigroup $(Q,\ast )$ associated with $G$. (Recall that a
\emph{quasigroup} is a set with binary operation $\circ$ such that if in the
equation $x\circ y=z$ any two of the elements are known then the third is
uniquely determined. For finite quasigroups this is equivalent to the
multiplication table being a latin square.) The quasigroup $Q$ has $G$ as its
underlying set, and the multiplication $\ast $ is given by $ g\ast h=g\cdot
h^{-1}$. Then $(Q,\ast )$ clearly satisfies the identity
\begin{equation}
(x\ast y)\ast (z\ast y)=x\ast z,  \label{Eq:Dedekind}
\end{equation}
and conversely if $Q$ is a quasigroup satisfying $(\ref{Eq:Dedekind})$ we can
reconstruct a group from $Q$. Quasigroups $(Q,*)$ satisfying
$(\ref{Eq:Dedekind})$ are known as \emph{Ward quasigroups}.

As mentioned above, ways of  axiomatizing a group  based on the right division
operation have appeared in several works, the first apparently being
\cite{Ward}. In Section \ref{Sc:WQ} we list some of these sets of axioms and
give concise proofs that several of the various identities which have appeared
are in fact equivalent. Our list is not exhaustive---for example \cite{PolTG}
contains several more equivalent identities. We also give some basic properties
of Ward quasigroups, some of which appear to be new.

We then discuss multiplication tables. We consider a finite group $G$ with
a cyclic subgroup $S$ of order $m$. If the elements of $G$ are ordered in a
special way using the left cosets of $S$, it is a consequence of
$(\ref{Eq:Dedekind})$ that the table which is obtained from the corresponding
Ward quasigroup $Q$ is a block matrix with $m\times m$ circulant blocks, in
which given any pair of rows the product of elements in the same column is
constant. This way of presenting a variation on the multiplication table of a
group appears to be new. For small groups this makes the description of the
multiplication table of $Q$ (and hence the group matrix) very concise. This
enables us to prove elementary facts about small groups combinatorially and we
also show that a group with a cyclic subgroup of index $2$ is completely
determined by the first row of the table. We call the permutation represented
by the first row of such a table an \emph{inverse pattern}, and we also show
that if $S$ is of index $3$, $Q$ and hence $G$ is determined by an inverse
pattern relative to $S$ and one other entry in the table.

More sophisticated tools have largely replaced multiplication tables of groups.
However when associativity is dropped many of these tools are no longer
available and often examples are constructed by indicating a multiplication
table. A generalization of Ward quasigroups is obtained when the operation
$x\ast y=x\cdot y^{-1}$ is based on a loop $(G,\cdot )$ with the
antiautomorphic inverse property. We indicate how symmetrical tables may be
constructed for all $6$ nonassociative Moufang loops of order at most $16$ with
the blocks being either circulants or reversed circulants. While we do not
pursue this here, there is an indication that new constructions of families of
Moufang loops can be obtained.

The implications for group matrices are discussed. We show how the block
circulant structure can be used to effect their partial diagonalization
combinatorially. This is equivalent to effecting a decomposition of the
regular representation into representations induced from $1$-dimensional
representations of $S$ up to $G$.

In the final section we indicate directions in which the work may lead.


\section{Ward Quasigroups}

\label{Sc:WQ}

\noindent Ward quasigroups have appeared in several different guises. The
concept (not the name) is due to M.~Ward \cite{Ward}. Rabinow \cite{Rab1}
discovered Ward quasigroups independently while axiomatizing groups using the
right division $x\cdot y^{-1}$ instead of $x\cdot y$. (Actually in \cite{Rab1}
he refers to a paper already submitted but we have been unable to discover a
published version.) The identity $(\ref{Eq:Dedekind})$ is mentioned for the
first time in Furstenberg \cite{Furstenberg}. The name \emph{Ward quasigroups}
was coined in 1978 by Cardoso and da Silva \cite{CardosoSilva}. Chatterjea
\cite{Chatterjea} and Polonijo \cite{Polonijo} were the first to show that Ward
quasigroups are exactly quasigroups satisfying $(\ref{Eq:Dedekind})$. The Ward
quasigroups corresponding to abelian groups, sometimes called \emph{subtractive
quasigroups}, are studied in \cite{Silva}, \cite{Morgado} and \cite{Whittaker}.

The following theorem generalizes \cite{Chatterjea} and provides the proof for
the claims made in the opening paragraphs of \cite{Rab1}. For additional
equivalent conditions, see \cite{PolTG}.

\begin{theorem}[Ward Quasigroups]
\label{Th:Ward} Let $G$ be a set and $*$ a binary operation on $G$. Then the
following conditions are equivalent:

\begin{enumerate}
\item[(i)] $G$ is a quasigroup satisfying $(1)$.

\item[(ii)] The operation $*$ satisfies $(1)$ and $a*G=G$ for every $a\in G$
$($cf.\ \emph{\cite{Furstenberg})}.

\item[(iii)] The square $a*a=e$ is independent of $a\in G$,
\begin{equation}  \label{Eq:OrigW}
(a*b)*c = a*(c*(e*b))
\end{equation}
holds for every $a$, $b$, $c\in G$, and if $e*a = e*b$ for some $a$, $b\in G$
then $a=b$ $($cf.\ \emph{\cite{Ward})}.

\item[(iv)] The square $a*a=e$ is independent of $a\in G$. Let $a^{\prime}=a*e$
for $a\in G$. Then $a^{\prime\prime}=a$ and
\begin{equation}  \label{Eq:OrigR}
(a*b^{\prime})*c^{\prime}=a*(b*c^{\prime})^{\prime}
\end{equation}
for every $a$, $b$, $c\in G$ $($cf.\ \emph{\cite{Rab1})}.

\item[(v)] Let $a\cdot b=a*((b*b)*b)$. Then $(G,\cdot)$ is a group. Its
neutral element $e=a*a$ is independent of $a$, the inverse of $a\in G$ is
given by $a^{-1}=e*a$, and $a\cdot b^{-1}=a*b$ holds for every $a$, $b\in G$.
\end{enumerate}
\end{theorem}

\begin{proof}
(i) implies (ii). All translations are bijections in a quasigroup.

(ii) implies (iii). Given $a$, $b\in G$, there is $c\in G$ such that $a*c=b$,
since $a*G=G$. Then, by $(\ref{Eq:Dedekind})$, $b*b=(a*c)*(a*c)=a*a$ is
independent of $b$, and we call it $e$.

Note that $a*e=a$ for every $a\in G$, since $a=a*b$ for some $b\in G$, and
thus $a*e=(a*b)*(b*b)=a*b=a$, by $(\ref{Eq:Dedekind})$.

Before we deduce $(\ref{Eq:OrigW})$, we show that $G*a=G$ for any $a\in G$. Let
$a$, $b\in G$. There is $c\in G$ such that $b=a*c$, and, in turn, there is
$d\in G$ such that $a=c*d$. By $(\ref{Eq:Dedekind})$,
$b=a*c=(a*d)*(c*d)=(a*d)*a$, and thus $b\in G*a$.

We need to show that $(a\ast b)\ast c=a\ast (c\ast (e\ast b))$. By the previous
paragraph, there exists $d\in G$ such that $c=d\ast b$. Then $ (a\ast b)\ast
c=(a\ast b)\ast (d\ast b)=a\ast d=a\ast (d\ast e)=a\ast ((d\ast b)\ast (e\ast
b))=a\ast (c\ast (e\ast b))$.

Finally assume that $e*a = e*b$ for some $a$, $b\in G$. Then $a*e =
(a*a)*(e*a) = (b*b)*(e*b) = b*e$, too. Since $a=a*e=b*e=b$, we are done.

(iii) implies (iv). We proceed similarly to \cite{Ward}. Note that $e=e*e$ by
the uniqueness of $e$. Hence $e*a=(e*e)*a = e*(a*(e*e)) = e*(a*e)$, where the
middle equality follows by $(\ref{Eq:OrigW})$. Consequently, $a=e*a$, since we
are allowed to cancel $e$ on the left. Then $e*a = (e*a)*e = e*(e*(e*a))$ by
$(\ref{Eq:OrigW})$, and therefore $a=e*(e*a)=a^{\prime\prime} $. Using
$(\ref{Eq:OrigW})$ and $a=a^{\prime\prime}$ repeatedly we obtain,
$(a*b^{\prime})*c^{\prime}= (a*(e*b))*(e*c) = a*((e*c)*(e*(e*b))) = a*((e*c)*b)
= a*(e*(b*(e*c))) = a*(b*c^{\prime})^{\prime}$, proving $(\ref {Eq:OrigR})$.

(iv) implies (v). Define $a\cdot b =a*((b*b)*b)=a*(e*b)=a*b^{\prime}$. By
$(\ref{Eq:OrigR})$, $(a\cdot b)\cdot c=(a*b^{\prime})*c^{\prime}=
a*(b*c^{\prime})^{\prime}=a\cdot (b\cdot c)$. We have $e\cdot
a=e*a^{\prime}=a^{\prime\prime}=a$, and, by $(\ref{Eq:OrigR})$ again, $a\cdot
e=a*e^{\prime}=a*(a^{\prime}*a^{\prime})^{\prime}=
(a*a^{\prime\prime})*a^{\prime}= (a*a)*a^{\prime}=
e*a^{\prime}=a^{\prime\prime}=a$. Furthermore, $a\cdot
a^{\prime}=a*a^{\prime\prime}=a*a=e$ and $a^{\prime}\cdot
a=a^{\prime}*a^{\prime}=e$. Therefore $(G,\cdot)$ is a group. As
$b=(b^{-1})^{-1}$ in any group, we have $a\cdot b^{-1} = a*(b^{-1})^{-1}=a*b$.

(v) implies (i). Since $a*b=a\cdot b^{-1}$, we have $(a*c)*(b*c)=(a\cdot
c^{-1})\cdot (b\cdot c^{-1})=a\cdot b^{-1}=a*b$. The equation $a*b=c$ can be
written as $a\cdot b^{-1}=c$, and therefore has a unique solution in $G$
anytime two of the three elements $a$, $b$, $c\in G$ are given.
\end{proof}

\begin{remark}
Equations $(\ref{Eq:OrigW})$ and $(\ref{Eq:OrigR})$ translate to the following
respective group identities, using (v):  $(a\cdot b^{-1})\cdot c^{-1}=a\cdot
(c\cdot (b^{-1})^{-1})^{-1}$ and $(a\cdot (b^{-1})^{-1})\cdot
(c^{-1})^{-1}=a\cdot ((b\cdot (c^{-1})^{-1})^{-1})^{-1}$. They are therefore
convoluted versions of the associative law and properties of ${}^{-1}$. Also
note that Furstenberg, Ward and Rabinow do not assume that the underlying
groupoid is a quasigroup. A groupoid satisfying $(\ref{Eq:Dedekind})$ is called
a \emph{T-groupoid} in \cite{PolTG}. The identity $(\ref{Eq:Dedekind})$ is
often called a \emph{right transitive identity}.
\end{remark}

The equations
\begin{equation*}
a\cdot b= a*((b*b)*b),\quad\quad a*b = a\cdot b^{-1}
\end{equation*}
of Theorem \ref{Th:Ward} show how to convert a Ward quasigroup to a group and
vice versa. Hence the essence of Ward quasigroups is the replacement of the
ordinary group multiplication $a\cdot b$ with the right division $a*b=a\cdot
b^{-1}$, as was observed already in \cite{Ward}, \cite{Rab1}. There is a Galois
correspondence between the two operations:

\begin{lemma}
Denote by $\mathrm{Wa}(G)$ the Ward quasigroup constructed from the group $G$,
and by $\mathrm{Gr}(Q)$ the group constructed from the Ward quasigroup $Q$.
Then $\mathrm{Gr}(\mathrm{Wa}(G))=G$ for every group $G$, and
$\mathrm{Wa}(\mathrm{Gr}(Q))=Q$ for every Ward quasigroup $Q$.
\end{lemma}

\begin{proof}
Let $*$ be the multiplication in a Ward quasigroup $Q$, $\cdot$ the
multiplication in $\mathrm{Gr}(Q)$, and $\circ$ the multiplication in
$\mathrm{Wa}(\mathrm{Gr}(Q))$. Then $x\circ y = x\cdot y^{-1} =
x*(y^{-1})^{-1}=x*y$. Similarly for $\mathrm{Gr}(\mathrm{Wa}(G))=G$.
\end{proof}

Ward quasigroups are therefore in one-to-one correspondence with groups, and
can be used to offer new insight into groups. From now on we will use the
term Ward quasigroup  to describe $\mathrm{Wa}(G)$, where $G$ is a group.
Multiplication in $G$ will be written as $ab$ instead of $a\cdot b$.

Following Rabinow's notation, when $(Q,*)$ is a Ward quasigroup with $e=a*a$,
let us define the bijection $^{\prime}:Q\to Q$ by $a\mapsto a^{\prime}=e*a$.
Note that $(a*b)^{\prime}=e*(a*b)=e(a*b)^{-1}=(ab^{-1})^{-1}=ba^{-1}=b*a$, and
$a^{\prime\prime}=a$, by Theorem \ref{Th:Ward}. Thus
$aa^{\prime}=a*a^{\prime\prime}=e$, and $a^{-1}=a^{\prime}$ follows.

We list some additional properties of Ward quasigroups.

\begin{lemma}
Let $(Q,*)$ be a Ward quasigroup. Then $\mathrm{Gr}(Q)$ is a commutative
group if and only if $\ ^{\prime}$ is an automorphism of $(Q,*)$.
Conversely, the Ward quasigroup $Q=\mathrm{Wa}(G)$ is commutative if and
only if $G$ is an elementary abelian $2$-group.
\end{lemma}

\begin{proof}
We have $(a*b)^{\prime}=ba^{-1}$ and
$a^{-1}b=a^{-1}*b^{-1}=a^{\prime}*b^{\prime}$. Thus $\mathrm{Gr}(Q)$ is
commutative if and only if $\ ^{\prime}$ is an automorphism of $(Q,*)$.

Conversely, $a*b^{\prime}= ab$ and $b^{\prime}*a=b^{-1}a^{-1}=(ab)^{-1}$
show that $\mathrm{Wa}(G)$ is commutative if and only if every element of $G$
is of exponent $2$.
\end{proof}

For any quasigroup $(Q,*)$, the \emph{associator} $[x,y,z]$ of $x$, $y$, $z\in
Q$ is the unique element $w$ such that $(x*(y*z))*w = (x*y)*z$.

\begin{lemma}
Let $(Q,*)$ be a Ward quasigroup. Then $[x,y,z]= z*((y*z)*y)$. In
particular, $[x,y,z]$ is independent of $x$.
\end{lemma}

\begin{proof}
Let $(Q,*)=W(G)$ be a Ward quasigroup, and $x$, $y$, $z\in Q$. If $w$ is
such that $(x*(y*z))*w=(x*y)*z$ then $x(yz^{-1})^{-1}w^{-1}=xy^{-1}z^{-1}$,
or $w=zyzy^{-1} = z(yz^{-1}y^{-1})^{-1} = z*((y*z)*y)$.
\end{proof}

The following consequence of $(\ref{Eq:Dedekind})$ was observed by
J.~D.~Phillips:

\begin{lemma}
Ward quasigroups satisfy the \emph{right semimedial law}:
\begin{equation}
(x*y)*(z*y)=(x*z)*(y*y).
\end{equation}
\end{lemma}

We conclude this section with a result concerning the identity
$(\ref{Eq:Dedekind})$ and generators of a quasigroup $Q$. The first part of
Lemma \ref{Lm:Gens} is due to Polonijo \cite{PolQQ}. He calls the elements of
$Y(Q)$ \emph{right quasiunits of $Q$}.

\begin{lemma}\label{Lm:Gens}
Let $Q=(Q,*)$ be a quasigroup $($not necessarily Ward$)$, and let $Y(Q)=\{y\in
Q;\;(x*y)*(z*y)=x*z$ for every $x$, $z\in Q\}$. If $Y(Q)$ is nonempty, it is a
subquasigroup of $Q$. Consequently, if $X$ is a generating subset of $Q$ such
that $X\subseteq Y(Q)$ then $Q$ is a Ward quasigroup.
\end{lemma}

\begin{proof}
Pick $y_1$, $y_2\in Y=Y(Q)$ and $x$, $z\in Q$. Then there are $x^{\prime}$,
$z^{\prime}\in Q$ such that $x=x^{\prime}*y_2$, $z=z^{\prime}*y_2$. Therefore
\begin{multline*}
(x*(y_1*y_2))*(z*(y_1*y_2)) =
((x^{\prime}*y_2)*(y_1*y_2))*((z^{\prime}*y_2)*(y_1*y_2)) \\
=(x^{\prime}*y_1)*(z^{\prime}*y_1) = x^{\prime}*z^{\prime}=
(x^{\prime}*y_2)*(z^{\prime}*y_2) = x*z,
\end{multline*}
and $Y$ is a subquasigroup. The rest follows.
\end{proof}


\section{Multiplication Tables}

\noindent In this section, we will restrict our attention to finite Ward
quasigroups.

Let $(Q,*)=\mathrm{Wa}(G)$ be a Ward quasigroup of order $n$, and let $S$ be a
cyclic subgroup of $G$ of order $m$ with generator $s$. Then $S$ is a
subquasigroup of $Q$ and the elements of $S$ can be listed as $e$, $s$, $s^2$,
$\dots$, $s^{m-1}$, where the powers are calculated in $G$.

Let $k=n/m$. Assume that $a_{1}=e$, $a_{2}$, $\dots $, $a_{k}$ form a set of
representatives of the left cosets $\{gS;\;g\in G\}$ of $S$ in $G$. Let us
construct a multiplication table $M$ of $Q$ as follows: order the elements of
the coset $a_{i}S$ as $a_{i}$, $a_{i}s$, $\dots $, $a_{i}s^{m-1}$. Then order
all elements of $Q$ by first using the elements of $a_{1}S$, then $a_{2}S$,
etc. This ordering will be used to label both rows and columns of $M $. (Thus
the set of elements in the $(i,j)$th block of the table is $a_{i}Sa_{j}^{-1}$.)

\begin{proposition}
\label{Pr:Table} Let $M$ be the multiplication table of $Q$ as described
above. Then

\begin{enumerate}
\item[(i)] $M=(m_{ij})$ consists of $k^2$ circulant matrices $C_{ij}$, each
of size $m$;

\item[(ii)] if we take any pair of rows of $M$, the product of each two
entries in the same column is constant, i.e., $m_{ij}*m_{kj} = m_{il}*m_{kl}$
for every $i$, $j$, $k$, $l$;

\item[(iii)] if the $j$th column of $M$ is labelled by $q\in Q$, then
$m_{1j}=q^{-1}$;

\item[(iv)] all the diagonal elements of $M$ are equal to $e$;

\item[(v)] the transpose of $C_{ij}$ is $(C_{ji})^{\prime }$. Here if
$A=(a_{i,j})$ is a matrix we use $A^{\prime }$ to denote the matrix
$(a_{i,j}^{\prime })$
\end{enumerate}
\end{proposition}

\begin{proof}
A circulant of order $m$ is determined by the following property: an entry
in the $(i,j)$th position is equal to the entry in the $(i+1,j+1)$th
position, where i+1 and j+1 are reduced modulo $m$. In the block $C_{ij}$,
if the $(k,l)$th entry is $x*y$ the $(k+1,l+1)$th entry is $(xs)*(ys) =
(x*s^{\prime})*(y*s^{\prime})=x*y$, where we again reduce $k+1$ and $l+1$
modulo $m$. Thus every block $C_{ij}$ is a circulant matrix, and we have
shown (i).

Assume that the $j$th column is labelled by $q$. Then $m_{ij}*m_{kj}
=(m_{i1}*q)*(m_{k1}*q)=m_{i1}*m_{k1}$, which shows (ii). Moreover, $m_{1j} =
e*q = q^{\prime}= q^{-1}$, which shows (iii). By Theorem \ref{Th:Ward}, $x*x=e$
for every $x\in Q$, and (iv) follows. Finally, $(x*y)=(y*x)^{\prime}$ implies
(v).
\end{proof}
\begin{remark}
If the table for a Ward quasigroup is constructed with any ordering of the
elements then condition (ii) is satisfied, and conversely if any quasigroup
table satisfies (ii) then the quasigroup is a Ward quasigroup. However with our
specific ordering described above to test the table for (ii) it is sufficient
to test only pairs of rows which correspond to the first line of any circulant
block, i.e. the rows in the $im^{th}$ places for $i=1,\dots,n/m$.
\end{remark}

\begin{example}
\label{Ex:S3} Let $G$ be the symmetric group on three elements, and let $S$
be the unique cyclic subgroup of order $3$ in $G$. Let $e=1$, $2$, $3$
denote the elements of $S$. Since every element of $G\setminus S$ is an
involution, Proposition \ref{Pr:Table} implies that the (incomplete)
multiplication table $M$ of $Q=\mathrm{Wa}(G)$ must be
\begin{equation*}
M=\begin{array}{c|ccc|ccc} \ast  & 1 & 2 & 3 & 4 & 5 & 6 \\ \hline
1 & 1 & 3 & 2 & 4 & 5 & 6 \\
2 & 2 & 1 & 3 & 6 & 4 & 5 \\
3 & 3 & 2 & 1 & 5 & 6 & 4 \\ \hline
4 & 4 & 6 & 5 & 1 &  &  \\
5 & 5 & 4 & 6 &  & 1 &  \\
6 & 6 & 5 & 4 &  &  & 1
\end{array}
\quad .
\end{equation*}
Furthermore, using condition (ii) of Proposition \ref{Pr:Table} for rows $1$
and $3$ we deduce that $4\ast 5=1\ast 3=2$, and the complete table is
determined as
\begin{equation*}
M=
\begin{array}{c|ccc|ccc}
\ast  & 1 & 2 & 3 & 4 & 5 & 6 \\ \hline
1 & 1 & 3 & 2 & 4 & 5 & 6 \\
2 & 2 & 1 & 3 & 6 & 4 & 5 \\
3 & 3 & 2 & 1 & 5 & 6 & 4 \\ \hline
4 & 4 & 6 & 5 & 1 & 2 & 3 \\
5 & 5 & 4 & 6 & 3 & 1 & 2 \\
6 & 6 & 5 & 4 & 2 & 3 & 1
\end{array}
\quad .
\end{equation*}
\end{example}


\section{Inverse Patterns}

\noindent Given a group $G$ of order $n$, a cyclic subgroup $S$ of order $m$, a
set of representatives $a_1=e$, $\dots$, $a_{n/m}$ of left cosets of $S$ in
$G$, and an order in which the cosets are listed, the permutation defined by
the first row of $M$ will be referred to as an \emph{inverse pattern} (cf.\
Proposition \ref{Pr:Table}(iii)).

Every inverse pattern $\iota$ is an involution such that $\iota(S)=S$. When $S$
is normal in $G$ then $\iota(aS)=a^{-1}S$ for every coset $aS$.

\begin{example}
\label{Ex:9} Let $h$ be the permutation $h=(1)(23)(47)(58)(69)$. Note that $h$
is an involution that stabilizes the block $\{1$, $2$, $3\}$ and interchanges
the blocks $\{4$, $5$, $6\}$, $\{7$, $8$, $9\}$. It therefore appears to be a
candidate for an inverse pattern, of a group with a normal cyclic subgroup of
order $3$. However, we claim that $h$ is not an inverse pattern of any Ward
quasigroup $Q$ with blocks of size $3$.

The permutation $h$ forces the following entries in the multiplication table
of $Q$:
\begin{equation*}
M=
\begin{array}{c||ccc|ccc|ccc}
* & 1 & 2 & 3 & 4 & 5 & 6 & 7 & 8 & 9 \\ \hline\hline
1 & 1 & 3 & 2 & 7 & 8 & 9 & 4 & 5 & 6 \\
2 & 2 & 1 & 3 & 9 & 7 & 8 & 6 & 4 & 5 \\
3 & 3 & 2 & 1 & 8 & 9 & 7 & 5 & 6 & 4 \\ \hline
4 & 4 & 6 & 5 & 1 & 2 & 3 &  &  &  \\
5 & 5 & 4 & 6 & 3 & 1 & 2 &  & A &  \\
6 & 6 & 5 & 4 & 2 & 3 & 1 &  &  &  \\ \hline
7 & 7 & 9 & 8 &  &  &  & 1 & 2 & 3 \\
8 & 8 & 7 & 9 &  & B &  & 3 & 1 & 2 \\
9 & 9 & 8 & 7 &  &  &  & 2 & 3 & 1
\end{array}
\quad.
\end{equation*}
The corresponding group $\mathrm{Gr}(Q)$ then satisfies $4\cdot 8 =
4*h(8)=4*5=2$ and $8\cdot 4 = 8*h(4)=8*7=3$, which contradicts the fact that
every group of order $9$ is commutative.

We can also give a purely combinatorial argument. There are three choices
for the two unspecified blocks $A$, $B$ of $M$. If we assume that $4*7=7$,
we can use rows $3$, $4$ of $M$ to deduce $5*7=3*4=8$, $6*7=2*4=9$. Since
both blocks are circulants and the transpose of $A$ is $h(B)$, we can fill
them up. Similarly when $4*7=8$ or $4*7=9$. The three completions of $M$ are
\begin{equation*}
\begin{array}{ccc|ccc}
&  &  & 7 & 9 & 8 \\
&  &  & 8 & 7 & 9 \\
&  &  & 9 & 8 & 7 \\ \hline
4 & 5 & 6 &  &  &  \\
6 & 4 & 5 &  &  &  \\
5 & 6 & 4 &  &  &
\end{array},\quad
\begin{array}{ccc|ccc}
&  &  & 8 & 7 & 9 \\
&  &  & 9 & 8 & 7 \\
&  &  & 7 & 9 & 8 \\ \hline
5 & 6 & 4 &  &  &  \\
4 & 5 & 6 &  &  &  \\
6 & 4 & 5 &  &  &
\end{array},\quad
\begin{array}{ccc|ccc}
&  &  & 9 & 8 & 7 \\
&  &  & 7 & 9 & 8 \\
&  &  & 8 & 7 & 9 \\ \hline
6 & 4 & 5 &  &  &  \\
5 & 6 & 4 &  &  &  \\
4 & 5 & 6 &  &  &
\end{array}.
\end{equation*}
In the first case, condition (ii) on rows $3$ and $4$ implies $4*8$=$3*4$. But
$4*8=9$ and $3*4=8$. In the second case, again considering rows $3$ and $4$,
$6*7 = 3*4$. But $6*7=7$ and again $3*4=8$. Using the same rows, condition (ii)
is also violated in the third case since it implies that $6*8=3*4$, and $6*8=7$
whereas $3*4=8$.
\end{example}

We have seen in Examples \ref{Ex:S3} and \ref{Ex:9} that an inverse
pattern can contain a large amount of information about the multiplication
table $M$ when the blocks of $M$ are relatively large. We will see later in this
section that, not surprisingly, the inverse pattern does not specify $M$ in
general. Nevertheless, when $[G:S]=2$, $M$ is determined:

\begin{lemma}
Let $G$ be a group and $S$ a cyclic subgroup of index $2$. Then an inverse
pattern of $G$ relative to $S$ specifies the multiplication table $M$ of
$Q=\mathrm{Wa}(G)$.
\end{lemma}

\begin{proof}
Assume that $M$ has been constructed as in Proposition \ref{Pr:Table} and that
the rows and columns of $M$ are labelled $1$, $\dots$, $n=2m$. $M$ consists of
four blocks $C_{ij}, i,j =1,2$. Given the inverse pattern, three of these
blocks are obviously determined, namely $C_{11}$, $C_{12}$ and $C_{21}$. Every
column of $C_{12}$ contains all entries $m+1$, $\dots$, $n$. Using (ii) we may
determine the product of any pair $i$, $j$ with $i$, $j \in \{m+1,...,2m\}$
since the product of the corresponding pair in the first column of $M$ is
already known, and thus the $C_{22}$ block is also determined.
\end{proof}

If $S$ is a normal cyclic subgroup of index $3$ the inverse pattern almost
specifies the group.

\begin{lemma}
\label{Lm:Index3} Assume that $S$ is a cyclic normal subgroup of index $3$ in
$G$. Then the multiplication table of $Q=\mathrm{Wa}(G)$ is specified by the
inverse pattern and by one entry in the $(2,3)$-block, or in the $(3,2)$
-block.
\end{lemma}

\begin{proof}
The blocks $C_{11}$, $C_{12}$, $C_{13}$, $C_{21}$ and $C_{31}$ are specified by
the inverse pattern and by the condition (v) of Proposition \ref{Pr:Table}.
Since the elements of the two left cosets different from $S$ are interchanged
by $\ ^{\prime}$, the diagonal blocks $C_{22}$ and $C_{33}$ are also specified.
Once any entry in $C_{23}$ or $C_{32}$ is known, both blocks can be filled, as
indicated in Example \ref{Ex:9}.
\end{proof}

The following lemma shows that for any two abelian groups of odd order $n$
with a cyclic subgroup $S$ of order $m$ the inverse patterns relative to $S$
can be made to coincide. Hence the inverse pattern is far from determining
the multiplication table.

\begin{lemma}
\label{Lm:InversePattern} Let $G$ be an abelian group of odd order $n$, and
let $S$ be a cyclic subgroup of $G$ of order $m$. The elements of $G$ can
then be ordered so that the inverse pattern of $Q=\mathrm{Wa}(G)$ is
\begin{equation}  \label{Eq:InversePattern}
1,\ m,\ \dots,\ 2\ |\ 2m+1,\ 3m,\ \dots,\ 2m+2\ |\ m+1,\ 2m,\ \dots,\ m+2\
|\ \dots
\end{equation}
\end{lemma}

\begin{proof}
Let $h$ be the map $x\mapsto x^{-1}$. Then $h(S)=S$ and $h(aS)\ne aS$ for
$aS\ne S$, otherwise the odd-order group $G/S$ contains an involution. We have
$h(h(aS))=aS$, and the coset $aS$ can therefore be coupled with $h(aS)$ . We
are free to choose a representative of each coset. Assume that if $a$ is the
representative of $aS$ then $a^{-1}$ is the representative of $h(aS)$ . Then
the elements of $aS$ are ordered as $a$, $as$, $as^2$, $\dots$, $ as^{m-1}$,
where $s$ is some fixed generator of $S$. The elements of $h(aS)$ then must be
ordered as $a^{-1}$, $a^{-1}s$, $\dots$, $a^{-1}s^{m-1}$. Since
$h(as^k)=a^{-1}s^{m-k}$, we are done.
\end{proof}

\begin{lemma}
Let $G$ be a group of odd order $n$, and let $S$ be a cyclic central
subgroup of $G$ of order $m$. The elements of $G$ can be ordered so that the
inverse pattern of $Q=\mathrm{Wa}(G)$ is $(\ref{Eq:InversePattern})$.
\end{lemma}

\begin{proof}
The proof of Lemma \ref{Lm:InversePattern} goes through word for word.
\end{proof}

The following lemma indicates how two inverse patterns of a group $G$ with
respect to a normal cyclic subgroup $S$ must be related.

\begin{lemma}
Let $G$ be a group, $S$ normal subgroup of $G$, and $\iota$ an inverse
pattern for $G$ with respect to $S$. Then any other inverse pattern for $G$
with respect to $S$ can be obtained from $\iota$ by changing the order in
which the cosets $aS$ are listed, and by applying a simultaneous cyclic
shift to each pair of cosets $aS$, $a^{-1}S$.
\end{lemma}

\begin{proof}
Since the left and right cosets of $S$ coincide, we can assume that they are
listed in the same order. We examine the (possibly equal) cosets $aS$,
$a^{-1}S$. Assume that $a$ is the representative of $aS$ and $b$ is the
representative of $a^{-1}S$ defined by $\iota$. Then there is a permutation
$\pi$ on $\{0$, $\dots$, $m-1\}$ such that $(as^{k})^{-1}=bs^{\pi(k)}$. Let
$c=as^k\in aS$ be another representative of $aS$. Then
$(cs^l)^{-1}=(as^{k+l})^{-1}=bs^{\pi(k+l)}$, where we calculate the exponents
modulo $m$.
\end{proof}

\begin{example}
Let $G=\langle a,b: a^7=b^3=e, b^{-1}ab=a^2\rangle$ be the Frobenius group of
order $21$. If we denote the unique cyclic subgroup of order $7$ in $G$ by $S$,
and choose $ba$, $b^2a$ as representatives of the remaining two left cosets, we
calculate that an inverse pattern with respect to $S$ is
\begin{equation*}
1\ 7\ 6\ 5\ 4\ 3\ 2\ |\ 15\ 18\ 21\ 17\ 20\ 16\ 19\ |\ 8\ 13\ 11\ 9\ 14\ 12\
10\ .
\end{equation*}
We use this pattern as the first row of the multiplication table $M$ of
$Q=\mathrm{Wa}(G)$. By Lemma \ref{Lm:Index3}, it suffices to specify one more
entry to complete $M$. By calculation $8*15=15$, and then the table is
determined as {\small
\begin{equation*}
\begin{array}{lll}
C(\ 1,\ 7,\ 6,\ 5,\ 4,\ 3,\ 2) & C(15,18,21,17,20,16,19) & C(\
8,13,11,9,14,12,10) \\
C(\ 8,14,13,12,11,10,9) & C(\ 1,\ 4,\ 7,\ 3,\ 6,\ 2,\ 5) &
C(15,20,18,16,21,19,17) \\
C(15,21,20,19,18,17,16) & C(\ 8,11,14,10,13,9,12) & C(\ 1,\ 6,\ 4,\ 2,\ 7,\ 5,\
3)
\end{array}.
\end{equation*}
}
\end{example}


\section{Generalized Ward Quasigroups Associated with Loops}

\noindent In this section, we briefly discuss the situation when we start with
a non-associative loop instead of a group, and we will see that under special
circumstances some of the symmetry of the multiplication table of the
corresponding quasigroup is retained.

A loop $G$ with neutral element $e$ has \emph{two-sided inverses} if for any
$x\in G$ there is $x^{-1}\in G$ such that $xx^{-1}=x^{-1}x=e$. A loop with
two-sided inverses has the \emph{antiautomorphic inverse property} if
$(xy)^{-1}=y^{-1}x^{-1}$. A loop is \emph{diassociative} if any two elements
generate a group. A diassociative loop clearly has the antiautomorphic inverse
property. A \emph{Moufang loop} is a loop satisfying the identity $x(y(xz)) =
((xy)x)z$. It is well known that Moufang loops are diassociative (cf.
\cite{Pflugfelder}).

Let $G$ be a loop with the antiautomorphic inverse property. In a similar
fashion to that above we may associate a ``generalized Ward quasigroup''
$(Q,\ast )$ to $G$ by $x\ast y=xy^{-1}$. In general, the left cosets of a
subloop $S$ need not partition $G$, and even if they do, the set
$(a_{i}S)(a_{j}S)^{-1}$ may contain more than $|S|$ elements. In order to avoid
these difficulties  we assume that $S$ is a normal cyclic subgroup of $G$. The
resulting multiplication table of $(Q,\ast )$ satisfies (iii), (iv) and (v) of
Proposition \ref{Pr:Table} (exercise), but if $G$ is not associative (ii)(but
not necessarily (i)) must fail. (When it is assumed that $G$ is only a loop,
and the multiplication is defined by $x\ast y = xy_\rho$, where $yy_\rho = e$,
then $m_{1j}=q_\rho$, condition (iv) holds, but (v) does not necessarily hold.)

Small Moufang loops were examined and give rise to tables which have much of
the symmetry of those for groups. If we consider the table for the smallest
non-associative Moufang loop $M_{12}$ of order 12 with respect to the unique
subgroup of order 3 we obtain the following.

Given symbols $c_{1}$, $\dots $, $c_{n}$, denote by $R(c_{1},\dots ,c_{n})$ the
\emph{reversed circulant matrix}
\begin{equation*}
\left(
\begin{array}{cccc}
c_{n} & c_{n-1} & \cdots  & c_{1} \\
c_{n-1} & c_{n-2} & \cdots  & c_{n} \\
\vdots  & \vdots  &  & \vdots  \\
c_{1} & c_{n} & \dots  & c_{2}.
\end{array}
\right)
\end{equation*}
Note that any Latin square of order $n\leq 3$ is a circulant or reversed
circulant. The table of $(Q,\ast )$ is
\begin{equation*}
\begin{array}{llll}
C(1,3,2) & C(4,5,6) & C(7,8,9) & C(10,11,12) \\
C(4,6,5) & C(1,2,3) & R(10,12,11) & R(7,9,8) \\
C(7,9,8) & R(10,12,11) & C(1,2,3) & R(4,6,5) \\
C(10,12,11) & R(7,9,8) & R(4,6,5) & C(1,2,3).
\end{array}
\end{equation*}
Note that the first row and column (of blocks) and the blocks on the
diagonal are determined by diassociativity. The symmetrical nature of the
remaining blocks is to be remarked. The table obviously violates condition
(i) of \ref{Pr:Table} and it easy to see that (ii) also fails.  Also note
that the inverse pattern of $(Q,\ast )$ is impossible for groups, as there
is no group of order $12$ with $9$ involutions.

There are 5 nonassociative Moufang loops of order $16$ and each of them
possesses a cyclic normal subloop of order $4$ (cf.\ \cite{Goodaire}). We have
checked that the multiplication tables of the associated quasigroups can be all
written in such a way that every $4\times 4$ block in the first row, first
column, and along the main diagonal is a circulant, while every other block is
a reversed circulant. It is probably not typical that larger Moufang loops
which are extensions of cyclic groups have tables of this type. For example a
Moufang loop $32$ with a cyclic normal subloop of order $8$ gives rise to off
diagonal blocks which are neither circulants nor reverse circulants. We remark
that Chein's construction $M_{2n}(G,2)$ \cite{Chein}, produces many of the
small Moufang loops, and circulants of reversed circulants may be present
because the dihedral group of order $2m$ and the generalised quaternion group
of order $2^{m}$ have ordinary multiplication tables which (with respect to
suitable ordering) consist of blocks which are either circulants or reverse
circulants. Many Moufang loops of small order contain dihedral or generalized
quaternion groups as subloops of index $2$ (cf.\ \cite{Chein},
\cite{Goodaire}).

We leave this section with an example:

\begin{example}
Let $(Q,\ast )$ be a quasigroup whose multiplication table is
\begin{equation*}
\begin{array}{cc}
C(1,3,2) & C(4,5,6) \\
C(4,6,5) & C(1,3,2)
\end{array}
\end{equation*}
Then the loop $G$ whose multiplication table is obtained from that above by permuting columns 2 and 3
 is the smallest nonassociative loop such that the
multiplication table for its generalized Ward quasigroup satisfies all the conditions of \ref{Pr:Table}
except condition (ii).
\end{example}


\section{The Group Matrix, Partial Diagonalization and Induced
Representations}

\noindent If $G$ is a finite group with associated Ward quasigroup $Q$, the
group matrix $X_G$ may be obtained from the multiplication table of $Q$ by
replacing each element $g$ by the variable $x_g$. From the results of Section 3
it follows that for every cyclic subgroup $S$ of $G$ with a compatible
ordering, $X_G$ is a block matrix $X_G=(B_{ij})_{r\times r}$, where each
$B_{ij}$ is a circulant of the form $C(x_{g_{k_1}}$, $\dots$, $ x_{g_{k_m}})$
where $g_{k_1}, \dots, g_{k_m}$ are elements of $G$. We denote this special
group matrix by $D_G(S)$, or $D_G$ if no ambiguity occurs.

\begin{example}
From Example \ref{Ex:S3}, the group matrix $D_{G}(C_3)$, $G=S_{3}$ is
\begin{equation*}
\left(
\begin{array}{ll}
C(x_{1},x_{3},x_{2}) & C(x_{4},x_{5},x_{6}) \\
C(x_{4},x_{6},x_{5}) & C(x_{1},x_{2},x_{3})
\end{array}
\right) .
\end{equation*}
\end{example}

\begin{lemma}
\label{Lm:GroupMatrix} For each circulant $C=C(a_1,\dots,a_m)$, if
\begin{equation*}
P=\left(
\begin{array}{ccccc}
1 & 1 & 1 & \dots & 1 \\
1 & \rho & \rho^2 & \dots & \rho^{m-1} \\
1 & \rho^2 & \rho^4 & \dots & \rho^{2(m-1)} \\
\vdots & \vdots & \vdots &  & \vdots \\
1 & \rho^{m-1} & \rho^{2(m-1)} & \dots & \rho^{(m-1)^2}
\end{array}
\right),
\end{equation*}
where $\rho=e^{2\pi i/n}$, then $P^{-1}CP$ is the diagonal matrix
\begin{equation*}
\mathrm{diag}(a_1+a_2+\cdots+a_m,\ a_1+\rho a_1+\cdots+\rho^{m-1}a_m,\dots,
a_1+\rho^{m-1}a_2+\cdots+\rho^{(m-1)^2}a_m).
\end{equation*}
\end{lemma}

\begin{proof}
This is by checking directly that each column of $P$ is an eigenvector of $
C(a_1$, $\dots$, $a_m)$.
\end{proof}

We denote by $\Lambda(a_1,\dots,a_m)$ the diagonal matrix $PC^{-1}P$ in
Lemma \ref{Lm:GroupMatrix}, and let $H$ be the $n\times n$ block diagonal
matrix $\mathrm{diag}(P,P,\dots,P)$ ($r$ blocks).

\begin{proposition}
If $G$ is any group of order $n=mr$ with cyclic subgroup $S$ of order $m$
then
\begin{equation*}
H^{-1}D_GH=\left(
\begin{array}{ccc}
\Lambda_{11} & \dots & \Lambda_{ir} \\
\vdots &  & \vdots \\
\Lambda_{r1} & \dots & \Lambda_{rr}
\end{array}
\right),
\end{equation*}
where $\Lambda_{ij}=\Lambda(x_{g_{{k_1}}}$, $\dots$, $x_{g_{k_m}})$ for all
$i$, $j$ and where $g_{k_1}, \dots, g_{k_m}$ are elements of $G$.
\end{proposition}

\begin{proof}
Each $m\times m$ block in the product $H^{-1}D_GH$ is $P^{-1}B_{ij}P$ and
the result follows directly.
\end{proof}

We continue the discussion of the example of $G=S_{3}$. Here
\begin{equation*}
H^{-1}D_GH=\left[
\begin{array}{cc}
\Lambda(x_{1},x_{3},x_{2}) & \Lambda(x_{4},x_{5},x_{6}) \\
\Lambda(x_{4},x_{5},x_{6}) & \Lambda(x_{1},x_{2},x_{3})
\end{array}
\right] .
\end{equation*}

If $\sigma$ is the permutation $(2,3,5,4)$ it is easily seen that permuting
the rows and columns of $H^{-1}DH$ by $\sigma$ is equivalent to conjugating
by a permutation matrix $R$, and the matrix we obtain is the block diagonal
matrix with blocks
\begin{eqnarray*}
D_{1}&=&\left[
\begin{array}{cc}
x_{1}+x_{2}+x_{3} & x_{4}+x_{5}+x_{6} \\
x_{4}+x_{5}+x_{6} & x_{1}+x_{2}+x_{3}
\end{array}
\right], \\
D_{2}&=&\left[
\begin{array}{cc}
x_{1}+\omega x_{2}+\omega^{2}x_{3} & x_{4}+\omega x_{5}+\omega^{2}x_{6} \\
x_{4}+\omega^{2}x_{5}+\omega x_{6} & x_{1}+\omega^{2}x_{2}+\omega x_{3}
\end{array}
\right], \\
D_{3}&=&\overline{D_{2}}.
\end{eqnarray*}

To any such block there is a representation of $G$. The explicit matrix
representing an element $g$ of $G$ is obtained by inserting $x_{g}=1$,
$x_{h}=0$ in the appropriate block. The representation is irreducible if and
only if the determinant of the block is irreducible as a polynomial in
$\{x_{1},...,x_{6}\}$. It is easily seen that det($D_{1}$)$
=u^{2}-v^{2}=(u+v)(u-v)$ where $u=x_{1}+x_{2}+x_{3}$, $v=x_{4}+x_{5}+x_{6}$,
which actually confirms that $D_{1}$ corresponds to the direct sum of the
trivial representation and the sign representation of $G$. Again by Frobenius'
theory since det$(D_{2})=$ det$(D_{3})$ the corresponding representations are
equivalent. We have rediscovered that there are three irreducible
representations of $G$, corresponding to the well-known character table
\begin{equation*}
\begin{array}{ccc}
1 & \{2,3\} & \{4,5,6\} \\ \hline
1 & 1 & 1 \\
1 & 1 & -1 \\
2 & -1 & 0
\end{array}
.
\end{equation*}

The above may be easily extended to any dihedral or generalized quaternion
group.

Consider again the Frobenius group $G$ of order $21$. Using the multiplication
table of the Ward quasigroup $Q=\mathrm{Wa}(G)$, we can see that $H^{-1}D_GH=K$
is obtained from the matrix $M$ described in Section 3 by replacing each
circulant $C(i_1$, $\dots$, $i_m)$ by $\Lambda(x_{i_1}$, $\dots$, $x_{i_m})$.

Let $\pi$ be the permutation
\begin{equation*}
\pi=(1)(2,4,10,8)(3,7,19,15)(5,13,17,9)(6,16)(11)(12,14,20,18)(21),
\end{equation*}
and $R_\pi$ the permutation matrix such that $(R_\pi)_{ij}$ equals $1$ if
$\pi(i)=j$, and $0$ otherwise. Then $R_\pi^{-1}KR_\pi$ is a block diagonal
matrix $\mathrm{diag}(B_1$, $\dots$, $B_7)$, where each $B_{s+1}$, $s=0$,\dots,$6$ is a $3\times
3$ block of the form {\footnotesize
\begin{equation*}
\left(
\begin{array}{lll}
\mu_s(1,7,6,5,4,3,2) & \mu_s(15,18,21,17,20,16,19) &
\mu_s(8,13,11,9,14,12,10) \\
\mu_s(8,14,13,12,11,10,9) & \mu_s(1,4,7,3,6,2,5) &
\mu_s(15,20,18,16,21,19,17) \\
\mu_s(15,21,20,19,18,17,16) & \mu_s(8,11,14,10,13,9,12) & \mu_s(1,6,4,2,7,5,3)
\end{array}
\right),
\end{equation*}
} where
\begin{equation*}
\mu_s(i_1,i_7,i_6,\dots,i_2) = x_{i_1}+\rho^sx_{i_2}
+\rho^{2s}x_{i_3}+\rho^{3s}x_{i_4}+\rho^{4s}x_{i_5}
+\rho^{5s}x_{i_6}+\rho^{6s}x_{i_7},
\end{equation*}
and $\rho=e^{(2 \pi i)/7}$. In fact, this decomposition enables us to
decompose the regular representation of $F_{21}$ into $7$ representations of
degree $3$, which each correspond to $B_{i}$, $i=1$, $\dots$, $7$. The
matrix representing the element $g$ is obtained by replacing $x_g$ by $1$
and $x_h$ by $0$, for $h\ne g$, in $B_i$. For instance, the matrices which
represent the generators $2$ and $8$ are
\begin{equation*}
\left(
\begin{array}{ccc}
\rho^s & 0 & 0 \\
0 & \rho^{2s} & 0 \\
0 & 0 & \rho^{4s}
\end{array}
\right)\quad\text{and}\quad \left(
\begin{array}{ccc}
0 & 0 & 1 \\
1 & 0 & 0 \\
0 & 1 & 0
\end{array}
\right),
\end{equation*}
respectively.

We use this information to obtain the character table of $F_{21}$. The block
$B_1$ corresponds to the representation
\begin{equation*}
2\mapsto I_3,\quad\quad 8\mapsto \left(
\begin{array}{ccc}
0 & 0 & 1 \\
1 & 0 & 0 \\
0 & 1 & 0
\end{array}
\right)=T,
\end{equation*}
and $T$ is similar to a diagonal matrix, the eigenvalues of $T$ being $1$,
$\omega$, $\omega^2$ (where $\omega=e^{2\pi i/3}$), and it splits into three
linear factors corresponding to the three linear representations of the group.
Each of the $6$ other blocks corresponds to one of the is two distinct
irreducible representations of degree $3$. The character table is
\begin{equation*}
\begin{array}{ccccc}
1 & \{2,3,5\} & \{4,6,7\} & 8-14 & 15-21 \\ \hline
1 & 1 & 1 & 1 & 1 \\
1 & 1 & 1 & \omega & \omega^2 \\
1 & 1 & 1 & \omega^2 & \omega \\
3 & \alpha & \beta & 0 & 0 \\
3 & \beta & \alpha & 0 & 0
\end{array}
.
\end{equation*}
where $\alpha=(-1+\sqrt{i})/2$ and $\beta=(-1-\sqrt{i})/2$. Analogously, for
an arbitrary group $G$, we obtain:

\begin{proposition}
If $G$ is an arbitrary group $G$ of order $n=mr$ with a cyclic subgroup $S$ of
order $m$, and $D_G$ is the group matrix corresponding to $Q_G$, then there
exists a permutation matrix $R_\pi$ such that $R_\pi^{-1}H^{-1}D_GHR_\pi$ is a
block diagonal matrix with $m$ blocks, each of size $r\times r$. This
effectively decomposes the regular representation of $G$ into a direct sum of
$m$ representations, each obtained by inducing an irreducible representation of
$S$ up to $G$.
\end{proposition}

We leave the proof to the reader. As a practical tool for group
representation theory the technique above would appear to be effective only
if the group has a cyclic subgroup of small index. Nevertheless, as a tool
to partially diagonalize the specializations of group matrices which occur
in Fourier analysis on finite groups or in the theory of tight frames, it
may have a wider application. We refer to \cite{dia} and \cite{wal} for
information on how group matrices appear in these contexts.


\section{Comments and questions}

\noindent Again we let $G$ be a group of order $n$ with a cyclic subgroup $S$
of order $m$.

1) It may be seen that the representation of $G$ which is induced from the
trivial representation of $S$ depends only on the sets of elements {$g_{k_1}$,
$\dots$, $g_{k_m}$} corresponding to the blocks $C(x_{g_{k_1}}$, $\dots$,
$x_{g_{k_m}})$. It may be interesting to relate these blocks to the theorem of
Artin on expressing any representation of $G$ in terms of representations
induced from cyclic subgroups.

2) Given a variety of groups the corresponding quasigroups must be
characterized by various identities. It may be interesting to examine these.

3) Is there a connection between coset enumeration with respect to $S$ and
our work here?

4) For an arbitrary group, is it possible to determine how much extra
information in addition to an inverse pattern with respect to a cyclic
subgroup determines the group?


\section{Acknowledgement}

\noindent We thank Michael Kinyon for providing us with references to the
early papers on Ward quasigroups.


\bibliographystyle{plain}

\begin{thebibliography}{99}

\footnotesize{

\bibitem{AT} L.~Auslander and P.~Tolmieri, \emph{Is computing
with finite Fourier transforms pure or applied mathematics?}, Bull. Amer. Math.
Soc. (N.S.) \textbf{1} (1979), 847--897.

\bibitem{CardosoSilva} J.~M.~Cardoso and C.~P.~da~Silva,
\emph{On Ward quasigroups}, An. \c Stiin\c t. Univ.\ ``Al. I. Cuza'' Ia\c si
Sec\c t. I a Mat. (N.S.) \textbf{24} (1978), no.\ 2, 231--233.

\bibitem{Chatterjea} S.~K.~Chatterjea, \emph{On Ward
quasigroups}, Pure Math.\ Manuscript \textbf{6} (1987), 31--34.

\bibitem[Chein(1978)]{Chein} O.~Chein, \emph{Moufang Loops of
Small Order}, Memoirs of the American Mathematical Society, Volume \textbf{13},
Issue 1, Number \textbf{197} (1978).

\bibitem{Silva} C.~P.~da~Silva, \emph{On a theorem of
Lagrange for Ward quasigroups} (Portuguese), Rev.\ Colombiana Ma.\ \textbf{12}
(1978), no.\ 3-4, 91--96.

\bibitem{dav} P.~J.~Davis, Circulant Matrices, Wiley, 1979.

\bibitem{dia} P.~Diaconis, \emph{Group Representations in
Probability and Statistics}, Lecture Notes, Institute of Mathematical
Statistics, Vol.\ \textbf{11}, Hayward, Ca.

\bibitem{F1} G.~Frobenius, \emph{\"Uber Gruppencharaktere},
Sitzungsber. Preuss. Akad. Wiss. Berlin (1896), 985-1021. (Gesammelte
Abhandlungen, (Springer-Verlag 1968), 1--37).

\bibitem{F2} G.~Frobenius, \emph{\"Uber die Primfactoren der
Gruppendeterminante}, Sitzungsber. Preuss Akad. Wiss Berlin (1896), 1343-1382.
(Gesammelte Abhandlungen, (Springer-Verlag 1968), 38--77).

\bibitem{Furstenberg} H.~Furstenberg, \emph{The inverse
operation in groups}, Proceedings of the American Mathematical Society
\textbf{6}, no.\ 6 (Dec., 1955), 991--997.

\bibitem[Goodaire, May and Raman(1999)]{Goodaire} Edgar
G.~Goodaire, Sean May, Maitreyi Raman, \emph{The Moufang Loops of Order less
than $64$}, Nova Science Publishers, 1999.

\bibitem{HJ1} H.-J.~Hoehnke and K.~W.~Johnson, \emph{The
3-characters are sufficient for the group determinant}, Proceedings of the
Second International Conference on Algebra, \emph{Contemporary Mathematics}
\textbf{184}(1995), 193-206.

\bibitem{J1} K.~W.~Johnson, \emph{On the group determinant},
Math. Proc. Cambridge Philos. Soc. \textbf{109} (1991), 299--311.

\bibitem{J2} K.~W.~Johnson, \emph{The Dedekind-Frobenius
group determinant, new life in an old method}, Proceedings, Groups St Andrews
97 in Bath, II, London Math. Soc. Lecture Note Series \textbf{261} (1999),
417-428.

\bibitem{Morgado} J.~Morgado, \emph{Defini\c{c}\~{a}o de quasigrupo
subtractivo por um \'unico axioma}, Gaz.\ Mat.\ (Lisboa)
\textbf{92}-\textbf{93} (1963), 17-18.

\bibitem{Pflugfelder} H.~O.~Pflugfelder, \emph{Quasigroups
and Loops: Introduction}, Heldermann Verlag, Berlin, 1990.

\bibitem{Polonijo} M.~Polonijo, \emph{A note on Ward
quasigroups}, An. \c Stiin\c t. Univ.\ ``Al. I. Cuza'' Ia\c si Sec\c t. I a
Mat. (N.S.) \textbf{32} (1986), no.\ 2, 5--10.

\bibitem{PolQQ} M.~Polonijo, \emph{Quasiunits in a quasigroup},
Rad Jugoslav. Akad. Znan. Umjet. No. \textbf{450} (1990), 117--121.

\bibitem{PolTG} M.~Polonijo, \emph{Transitive groupoids}, Portugaliae
Mathematica \textbf{50} (1993), no.\ 1, 63--74.

\bibitem{Rab1} D.~G.~Rabinow, \emph{Independent sets of
postulates for abelian groups and fields in terms of the inverse operations},
American Journal of Mathematics \textbf{59}, no. 1 (Jan., 1937), 211-224.

\bibitem{wal} R.~Reams, S.~Waldron, \emph{Isometric tight
frames}, Electron. J. Linear Algebra \textbf{9} (2002), 122--128

\bibitem{Vojt} P.~Vojt\v{e}chovsk\'y, \emph{The smallest
Moufang loop revisited}, Results in Mathematics \textbf{44} (2003), 189--193.

\bibitem{Ward} M.~Ward, \emph{Postulates for the inverse
operations in a group}, Transactions of the American Mathematical Society
\textbf{32}, no.\ 3 (Jul.,\ 1930), 520--526.

\bibitem{Whittaker} J.~V.~Whittaker, \emph{On the postulates defining a group},
Amer.\ Math.\ Monthly \textbf{62} (1955), 636--640.

}
\end{thebibliography}

\end{document}